\documentclass[reqno,11pt]{amsart}
\usepackage{amsmath,amssymb,latexsym, xcolor}

\def\pr{\text{ P\/}}
\def\ex{\text{E\/}}

\def\eps{\varepsilon}

\def\part{\partial}

\newtheorem{Theorem}{Theorem}[section]

\newtheorem{Corollary}[Theorem]{Corollary}

\theoremstyle{remark}

\numberwithin{equation}{section}

\begin{document}

\title[Connectivity of a directed core]{Counting strongly connected $(k_1,k_2)$-directed cores }
\date{\today}
\author{Boris Pittel}
\address{Department of Mathematics, The Ohio State University, Columbus, OH 43210-1175 (USA)}
\email{bgp@math.osu.edu}



\maketitle

\begin{abstract} Consider the set of all digraphs on $[N]$ with $M$ edges, whose minimum in-degree and minimum out-degree are at least $k_1$ and $k_2$ respectively. For $k:=\min\{k_1,k_2\}\ge 2$ and $M/N\ge \max\{k_1,k_2\}+\eps$, $M=\Theta(N)$, we show that, among those
digraphs, the fraction of $k$-strongly connected digraphs is $1-O\bigl(N^{-(k-1)})$. Earlier with Dan
Poole we identified a sharp edge-density threshold $c^*(k_1,k_2)$ for birth of a giant 
$(k_1,k_2)$-core in the random digraph $D(n,m=[cn])$. Combining the claims, for $c>c^*(k_1,k_2)$ with probability  $1-O\bigl(N^{-(k-1)})$ the giant $(k_1,k_2)$-core exists and is $k$-strongly connected.
\end{abstract}


\section{Results and related work}

Let the fixed integers $k_1\ge 2$, $k_2\ge 2$ be given. A digraph 
is called a directed $(k_1,k_2)$-core (dicore)
if its minimum in-degree and minimum out-degree are, at least, $k_1$ and $k_2$ respectively.

\begin{Theorem}\label{strcon}  Let $M$, $N$ be such that  $M\ge (\max \{k_1,k_2\}+\eps)N$, and
$M=O(N)$. 
For $\bold k:=(k_1,k_2)$, let $D_{\bold k}(N,M)$ be a dicore chosen uniformly at random among all simple dicores with $N$ vertices and $M$ directed eges. Then,  denoting $k=\min\{k_1,k_2\}$,   
$D_{\bold k}(N,M)$ is strongly connected with probability $1-O(N^{-2(k-1)})$, and $k$-strongly connected with probability $1-O(N^{-(k-1)})$.
\end{Theorem}

This theorem leads to a sharp asymptotic formula for $C_{\bold k}(N,M)$, the total number
of $k$-strongly connected dicores on $[N]$ with $M$ edges. Given $z>0$, let $\text{Poi}(z|k)$
stand for the $\text{Poisson}(z)$ conditioned on $\{\text{Poisson}(z)\ge k\}$. Introduce 
\[
p_k(z)=\pr\bigl(\text{Poisson}(z)\ge k\bigr),\quad f_k(z)=e^zp_k(z)=\sum_{j\ge k}\frac{z^j}{j!}.
\]
\begin{Corollary} \label{countcores} Under the conditions of Theorem \ref{strcon},
\begin{align*}
C_{\bold k}(N,M)&=\bigl(1+O(N^{-1}\log^6N)\bigr)\,\!\exp\!\left(\!-\frac{M}{N}-
\frac{1}{2}\prod_{j=1}^2\ex\bigl[\text{Poi}(z_j|k_j-1)\bigr]\!\right)\\
&\quad\times M!\prod_{j=1}^2\frac{f_{k_j}(z_j)^N}{z_j^M\sqrt{2\pi N\text{Var}(\text{Poi}(z_j|k_j))}},
\end{align*}
where $z_1$, $z_2$ satisfy $\ex\bigl[\text{Poi}(z_j|k_j)\bigr]=\frac{M}{N}$.
\end{Corollary}

\noindent The RHS expression is also an asymptotic formula for $D_{\bold k}(N,M)$ the total number of all
$(k_1,k_2)$-dicores, whose proof is a carbon copy of the formula (2.10) in Pittel \cite{Pit1} for
the special case $k_1=k_2=1$. (Without convergence rate, but under less restrictive condition
on $M=M(N)$, that formula for $D_{\bold k}(N,M)$ had been proved by P\'erez-Gim\'enez and Wormald \cite{PerWor}.)
So Corollary \ref{countcores} follows from
\[
 C_{\bold k}(N,M)/D_{\bold k}(N,M)=1-O\bigl(N^{-(k-1)}\bigr),
\]
a rephrased version of Theorem \ref{strcon}. To compare, it was proved in \cite{Pit1} that, for
$M-N\gg N^{2/3}$, $M=O(N)$,
\begin{multline*}
\mathcal C_{1,1}(N,M)=\bigl(1+O(N^2/(M-N)^3+N^{-1}\log^6N)\bigr)\,\!\exp\!\left(\!-\frac{M}{N}-
\frac{z^2}{2}\right)\\
\times \frac{M!}{2\pi N \text{Var}(\text{Poi}(z|1))}
\frac{f_1(z)^{2N}}{z^{2M}}\cdot\frac{\left(1-\tfrac{z}{f_1(z)}\right)^2}{1-\tfrac{z}{e^zf_1(z)}}\,
\exp\left[\frac{z}{f_1(z)}(2-e^{-z})\right],
\end{multline*}
where $z$ satisfies $\ex\bigl[\text{Poi}(z|1)\bigr]=\tfrac{M}{N}$. See 
\cite{PerWor} for a version of the last formula,  under broader conditions $M-N\to\infty$, $M=O(N\log N)$ but without convergence rates.  

It has long been known that, for the uniformly random digraph $D(n,m)$ with $n$ vertices and
$m$ edges $[n]$, the edge density $m/n=1$ is the sharp threshold for birth of a giant strong component, see Karp \cite{Kar} and T. \L uczak \cite{Luc2}, and a more recent  paper by T. \L uczak and 
Seierstad \cite{LucSei}. And in Pittel and Poole \cite{PitPoo1} it was proved that for the postcritical stage 
$m/n>1$ the number of vertices and the number of edges in the strong giant have a joint Gaussian distribution in the limit $n\to\infty$. This can be viewed as a directed analogue of the earlier 
result, Pittel and Wormald \cite{PitWor}, on the limit Gaussian distribution of the number of vertices
and the number of edges in the $2$-core, and the number of vertices in the ``forrest'' mantle,
of the giant component in the postcritical Erd\H os-R\'enyi random graph $G(n,m=[cn])$. 
Pittel, Spencer and Wormald \cite{PitSpeWor} determined the edge density threshold $c^*(k)$ for birth of a giant $k$-core ($k\ge 3)$ in $G(n,m=[cn])$:
\[
c^*(k)=\min_{z>0}\frac{z}{p_{k-1}(z)}.
\] 
Recently an analogous result for $D(n,m=[cn])$ was proved in Pittel and Poole \cite{PitPoo2}.  Let 
\[
c^*=c^*(\bold k):=
\min\limits_{z_1,z_2>0}\max\left\{\frac{z_1}{p_{k_1}(z_1)p_{k_2-1}(z_2)};
\frac{z_2}{p_{k_1-1}(z_1)p_{k_2}(z_2)}\right\}.
\]

\begin{Theorem}\label{birthcore} (\cite{PitPoo2}) Let $k_1,\,k_2\ge 0$, $\max\{k_1,k_2\}\ge 2$. Then
\begin{itemize}
\item for $c<c^* $, quite surely (q.s.) the $(k_1, k_2)$-core of $D(n, m=[c n])$ is empty;
\item for $c>c^* $, q.s. the $(k_1, k_2)$-core of $D(n, m=[c n])$ is not empty; in fact, there are
some $\alpha(c)=\alpha(\bold k,c)$ and $\beta(c) = \beta(\bold k, c)$, with $\alpha(c)>
\max\{k_1,k_2\}\beta (c)$, such that q.s. the $(k_1, k_2)$-core has $\beta(c) n + O(n^{1/2} \log n)$ vertices and $\alpha(c)n +O(n^{1/2}\log n)$ edges.
\end{itemize}
\end{Theorem}
Here ``quite surely'' means that the event in question has probability $1-O(n^{-K})$, for all $K>0$.
Our study left open an issue of strong connectivity of the giant $(k_1,k_2)$-core for $c>c^*$. Many
years ago T. \L uczak \cite{Luc1} proved that, for $k\ge 3$, w.h.p. if a $k$-core is present
in $G(n,m)$ it must be $k$-connected. Dan Poole has conjectured that likewise, for 
$k=\min\{k_1,k_2\}\ge 2$ and $c>c^*(\bold k)$,  w.h.p. the giant $(k_1,k_2)$-core of $D(n,m=[cn])$ is $k$-strongly
connected. The theorem \ref{strcon} can be used to confirm Dan's conjecture. Here is how.

The proof of Theorem \ref{birthcore} in \cite{PitPoo2} was based on analysis of a deletion algorithm: at each step a uniformly random vertex, with either light in-degree, i.e. below $k_1$, or light out-degree, i.e. below $k_2$, is deleted,
together with all incident edges. Instead of $D(n,m=[cn])$, we considered a uniformly random multi-digraph $\mathcal D(n,m=[cn])$, multiple loops and multiple edges allowed, on $n$ vertices with $m$ {\it labeled\/} edges. Conditioned on being simple, $\mathcal D(n,m=[cn])$ is distributed as $D(n,m=[cn])$.
For $m=O(n)$, $\mathcal D(n,m=[cn])$ is simple with positive asymptotic probability $e^{-c-c^2/2}$. Thus an event unlikely for $\mathcal D(n,m=[cn])$ is {\it equally\/} unlikely for $D(n,m=[cn])$. It is convenient
to view $\mathcal D(n,m=[cn])$ as a directed version of a sequence model invented by Chv\'atal \cite{Chv}
for study of $3$-colorability of $G(n,m)$, $m=O(n)$, and later used by Aronson, Frieze and Pittel \cite{AroFriPit} for analysis of a vertex deletion process at the heart of the Karp-Sipser greedy matching algorithm \cite{KarSip}.

Let us reproduce the definition of the sequence model  from \cite{PitPoo2}, since it
will be instrumental in our proofs in this paper as well.  Given a sequence $\bold x=(x_1,\dots,x_{2m})$, $x_i\in [n]$, we define a {\it multi\/}-digraph
$\mathcal D_{\bold x}$ with vertex set $[n]$ and (directed) edge set $\bigl[\{x_{2r-1},x_{2r}\}:\,1\le r\le m\bigr]$;
thus $e_{\bold x}(i,j)$, the number of directed edges $i\to j$, is $|\{r: x_{2r-1}=i,x_{2r}=j\}$. 
The in-degree sequence 
$\boldsymbol\delta_{\bold x}$
and the out-degree sequence $\boldsymbol\Delta_{\bold x}$ of $\mathcal D_{\bold x}$  are given by
$\delta_{\bold x}(i)=|\{r: x_{2r}=i\}|$, $\Delta_{\bold x}(i)=|\{r:x_{2r-1}=i\}|$, so that 
$\sum_{I\in [n]}\delta_{\bold x}(i)= \sum_{I\in [n]}\Delta_{\bold x}(i)=m$.
If $\bold x$ is distributed uniformly on the set $[n]^{2m}$ then $\mathcal D_{\bold x}$ and $\mathcal D(n,m)$ are
equi-distributed. Consequently $\boldsymbol\delta_{\bold x}$ and $\boldsymbol\Delta_{\bold x}$  are mutually independent, each distributed multinomially, with  $m$ trials and $n$ equally likely outcomes in each trial. 

The deletion algorithm delivers a sequence
$\{\bold x(t)\}$ where $\bold x(0)=\bold x$, and each $\bold x(t)\in ([n]\cup \{\star\})^{2m}$, where for all $r$, $x_{2r-1}(t)=\star$ if and only if $x_{2r}(t)=\star$. The $(\star,\star)$ pairs mark the 
locations $(2r-1,2r)$ in the original $\bold x(0)$ whose vertex occupants have been deleted after
$t$ steps. The process $\{\bold x(t)\}$ is obviously Markov, though the complexity of its sample
space makes it intractable. Let $\bold s_{\bold x}$ be a $[(k_1+1)(k_2+1)+1]$-tuple whose
components are the counts of vertices that are in/out-light, in-light/out-heavy, in-heavy/out-light,
in-heavy/out-heavy, and the total count of all edges in $\bold x$. We need that
many components since, to preserve Markovian property, we have to classify the in-light degrees and the out-light degrees according to their possible $k_1$ and $k_2$ values. Fortunately no similar classification is needed for the in-heavy degrees and the out-heavy degrees. It was proved in \cite{PitPoo2} that
the process $\{\bold s(t)\}:=\{s_{\bold x(t)}\}$ is indeed Markov, and that, conditioned on $\bold s(t)$,
the sequence $\bold x(t)$ is uniform. 

The upshot of this discussion is that, conditioned on the terminal vertex set and the terminal number of edges, the terminal sequence $\bold x$ is distributed uniformly.  So 
Theorem \ref{strcon} in combination  with  Theorem \ref{birthcore} from \cite{PitPoo2} yield
\begin{Corollary}\label{birthstr} Let $k_1\ge 2$, $k_2\ge 2$. If $c>c^*(\bold k)$ then with probability
$1-O\bigl(n^{-\min\{k_1,k_2\}+1}\bigr)$ 
the random digraph $D(n,m=[cn])$ has a giant $(k_1,k_2)$-core which is $\min\{k_1,k_2\}$-strongly connected.
\end{Corollary}

It is edge sparseness of the near-postcritical $D(n,m)$, i.e. $m$ being of order $n$, that forces us to push the in/out degrees
minimally upward from $1$. P\'erez-Gim\'enez and Wormald \cite{PerWor} proved that, when $m/n\to\infty$, the random digraph, whose all in/out degrees are merely positive, is strongly connected with high probability. 

Cooper and Frieze \cite{CooFri} studied a random directed graph with a {\it given\/} degree sequence,
which is a counterpart of the random undirected graph first introduced and analyzed by Molloy 
and Reed \cite{MolRee1}, \cite{MolRee2}. Among other results, it was proven in \cite{CooFri} that for a ``proper''
in/out {\it positive \/} degree sequence the random digraph has a giant strongly-connected
component comprised of almost all vertices.

\bigskip
\section{Proof of Theorem \ref{strcon}}
To analyze strong connectedness of $D_{\bold k}(N,M)$, we use  the Chv\'atal-type sequence model. Each admissible sequence $\bold x$ is obtained by filling the $M$ pairs of consecutive locations $(2r-1,2r)$, $r\in [M]$, 
with the vertices from $[N]$ such that every vertex appears at least $k_2$ times in the odd-numbered
locations and at least $k_1$ times in the even-numbered  locations. Let $\boldsymbol\delta$,
$\boldsymbol\Delta$ denote the in/out vertex degrees of an admissible sequence. Then
\begin{equation}\label{bdelta,bDelta=}
\delta_i\ge k_1,\,\,\Delta_i\ge k_2;\quad \sum_{i\in [N]}\delta_i=\sum_{i\in [N]}\Delta_i=M.
\end{equation}
As in \cite{PitPoo2}, the total number $S_{N,M}$ of these sequences is given by
\begin{equation}\label{mathcal SNM=}
\begin{aligned}
S_{N,M}&=\sum\limits_{\boldsymbol\delta,\boldsymbol\Delta\text{ meet }\eqref{bdelta,bDelta=}}\frac{(M!)^2}{\prod\limits_{i\in [N]}\delta_i!\,\Delta_i!}\\
&=(M!)^2 [z_i^Mz_o^M] f_{k_1}(z_i)^N\,f_{k_2}(z_o)^N\\
&=\Omega\left(N^{-1}(M!)^2\frac{f_{k_1}(z_i)^N\,f_{k_2}(z_o)^N}{z_i^M\,z_o^M}\right);
\end{aligned}
\end{equation}
here $z_i$ and $z_o$ are the roots of
\begin{equation}\label{zi,zo;roots}
\frac{z_i\,f_{k_1-1}(z_i)}{f_{k_1}(z_i)}= \frac{M}{N},\quad  \frac{z_o f_{k_2-1}(z_o)}{f_{k_2}(z_o)}=
\frac{M}{N}.
\end{equation}
The subindices ``i'' and ``o'' stand for ``in'' and  ``out'' respectively. The conditions \eqref{zi,zo;roots}
mean that $z_i$ and $z_o$ are the absolute minimum points for the functions
$z^{-M} f_{k_1}(z)^N$ and $z^{-M}f_{k_2}(z)^N$ respectively

We need to show that w.h.p. the uniformly random sequence is such that the corresponding
digraph is $k:=\min\{k_1,k_2\}$-strongly connected. To this end, we have to show that the  number of sequences $\bold x$, such that deletion of a set $T$ of $t<k$ vertices results in partition
of $[N]\setminus T$ into the disjoint union of a source/sink set  $A$ and a sink/source set $B=
([N]\setminus T)\setminus A$, is $o(S_{N,M})$ as $N\to\infty$. It suffices to consider 
the case when $|A|\le (N-t)/2$. Since we are interested in the sequences that can be induced by
a {\it simple\/} $(k_1,k_2)$-core, thus with each vertex of in-degree $k_1$ and out-degree $k_2$ at least, we may and will focus on $|A|\ge 2$.

For certainty, let $A$ be a source set.  Let us consider $t=0$ first.
Let $\nu\in [2,N/2]$, $\boldsymbol\mu=(\mu_1,\mu_{1,2},\mu_2)$,
$\mu_1+\mu_{1,2}+\mu_2=M$. Let  $S_{\nu,\boldsymbol\mu}$ be the total number of the
sequences $x$ that contain a source set $A$ of cardinality $\nu$,  such that $\mu_1$ ($\mu_2$ resp.) is the number of edges between the vertices in $A$ (in $B$ resp.), and $\mu_{1,2}$ is the number of edges from vertices in $A$ to
vertices in $B$.  It is necessary, of course,  that 
\begin{equation}\label{mularge}
\mu_1\ge k_1\nu,\,\,\mu_1+\mu_{1,2}\ge k_2\nu,\,\, \mu_{1,2}+\mu_2\ge k_1(N-\nu),
\,\, \mu_2\ge k_2(N-\nu).
\end{equation}
To be sure, the sequences $\bold x$ that contain several such source sets $A$ will be counted
more than once.

By symmetry,
\begin{equation}\label{Snu,mu=}
S_{\nu,\boldsymbol\mu}=\binom{N}{\nu}\binom{M}{\boldsymbol\mu}
\mathcal S_{\nu,\boldsymbol\mu};
\end{equation}
here $\mathcal S_{\nu,\boldsymbol\mu}$ is the total number of special admissible sequences 
with parameters $\nu$ and $\boldsymbol\mu$.  For an admissible sequence to be special,
the first  $\mu_1$ ordered
pairs must be formed by vertices from $[\nu]$, the next $\mu_{1,2}$ pairs--by pairs of vertices, left from $[\nu]$, right from $[N-\nu]:=
(\nu+1,\dots,N)$, and the last $\mu_2$ pairs--by vertices from $[N-\nu]$.  Further,
for the first block and $i\in [\nu]$, let $\delta_i^{(1)}$, ($\Delta_i^{(1)}$ resp.)
denote the number of pairs containing $i$ in the right slot (the left slot resp.); for the second block   and for $i\in [\nu]$, let $\Delta_i ^{(1,2)}$ denote the number of the left slots
containing $i$, and for $i\in [N-\nu]$, let $\delta_i^{(1,2)}$,  denote the number of right slots containing $i$; for the third block  and $i\in [N-\nu]$, let $\delta_i^{(2)}$, ($\Delta_i^{(2)}$ resp.) denote the number of pairs containing $i$ in the right slot (the left slot resp.). Then necessarily
\begin{equation}\label{nu,N-nu}
\begin{aligned}
\delta_i^{(1)}\ge k_1, \,\, \Delta_i^{(1)}+\Delta_i^{(1,2)}\ge k_2,\quad i\in [\nu], \\
 \delta_i^{(1,2)}+\delta_i^{(2)}\ge k_1,\,\, \Delta_i^{(2)}\ge k_2,\quad i\in [N-\nu],
 \end{aligned}
 \end{equation}
and
\begin{equation}\label{sums}
\begin{aligned}
\sum_{i\in [\nu]}\delta_i^{(1)}&=\sum_{i\in [\nu]}\Delta_i^{(1)}=\mu_1,\\
\sum_{i\in [N-\nu]}\!\delta_i^{(1,2)}&=\sum_{i\in [\nu]}\Delta_i^{(1,2)}=\mu_{1,2},\\
\sum_{i\in [N-\nu]}\!\delta_i^{(2)}&=\sum_{i\in [N-\nu]}\!\Delta_i^{(2)}=\mu_2.
\end{aligned}
\end{equation}
Enter the generating functions.  Introduce the indeterminates $\bold x=(x_1,x_{1,2},x_2)$, $\bold y=(y_1,y_{1,2},y_2)$,
and the notations $\boldsymbol\mu!=\mu_1!\mu_{1,2}!\mu_2!$, $\bold x^{\boldsymbol\mu}=
x_1^{\mu_1}x_{1,2}^{\mu_{1,2}}x_2^{\mu_2}$, $\bold y^{\boldsymbol\mu}= y_1^{\mu_1}y_{1,2}^{\mu_{1,2}}y_2^{\mu_2}$, and $x=x_{1,2}+x_2$, $y=y_1+y_{1,2}$. Using
\[
\sum_{d_1+\cdots+d_r\ge k}\,\prod_{j=1}^r\frac{\xi_j^{d_j}}{d_j!}=f_k\!\!\left(\sum_{j=1}^r\xi_j\right)\!\!,
\]
we write
\begin{multline*}
S_{\nu,\boldsymbol\mu}=\sum_{(\boldsymbol\delta,\boldsymbol\Delta)\text{ meet }\atop
\eqref{nu,N-nu}, \eqref{sums}}\frac{(\mu_1!)^2}{\prod\limits_{i\in [\nu]}(\delta_i^{(1)})!\cdot(\Delta_i^{(1)})!}\\
\times\frac{(\mu_{1,2}!)^2}{\prod\limits_{i\in [\nu],\, j\in [N-\nu]}(\Delta_i^{(1,2)})!\cdot(\delta_j^{(1,2)})!}\cdot
\frac{(\mu_2!)^2}{\prod\limits_{i\in [N-\nu]}(\delta_i^{(2)})!\cdot(\Delta_i^{(2)})!}\\
=(\boldsymbol\mu!)^2 \bigl[\bold x^{\boldsymbol\mu}\bold y^{\boldsymbol\mu}\bigr]
\left(\sum_{\delta\ge k_1}\frac{x_1^{\,\delta}}{\delta!}\right)^{\nu}\cdot\left(\sum_{\Delta\ge k_2}
\frac{y_2^{\,\Delta}}{\Delta!}\right)^{N-\nu}\\
\times\left(\sum_{\Delta^{(1)}+\Delta^{(2)}\ge k_2}\frac{y_1^{\,\Delta^{(1)}}y_{1,2}^{\,\Delta^{(1,2)}}}
{\Delta^{(1)}!\,\Delta^{(1,2)}!}\right)^{\nu}\cdot
\left(\sum_{\delta^{(1,2)}+\delta^{(2)}\ge k_1}\frac{x_{1,2}^{\,\delta^{(1,2)}}x_2^{\,\delta^{(2)}}}
{\delta^{(1,2}!\delta^{(2)}!}\right)^{N-\nu}\\
=(\boldsymbol\mu!)^2 
[x_1^{\mu_1}] f_{k_1}(x_1)^{\nu}\,[y_2^{\mu_2}] f_{k_2}(y_2)^{N-\nu}\\
\times [y_1^{\mu_1}y_{1,2}^{\mu_{1,2}}] f_{k_2}(y)^{\nu}\, [x_{1,2}^{\mu_{1,2}}x_2^{\mu_2}]f_{k_1}(x)^{N-\nu};
\end{multline*}
here
\begin{align*}
[y_1^{\mu_1}y_{1,2}^{\mu_{1,2}}] f_{k_2}(y)^{\nu}&=\binom{\mu_1+\mu_{1,2}}{\mu_1}
[y^{\mu_1+\mu_{1,2}}]f_{k_2}(y)^{\nu},\\
[x_{1,2}^{\mu_{1,2}}x_2^{\mu_2}]f_{k_1}(x)^{N-\nu}&=\binom{\mu_{1,2}+\mu_2}{\mu_2}
[x^{\mu_{1,2}+\mu_2}]f_{k_1}(x)^{N-\nu}.
\end{align*}
Thus
\begin{multline}\label{long}
\mathcal S_{\nu,\boldsymbol\mu}=\mu_1!\mu_2!(\mu_1+\mu_{1,2})!(\mu_{1,2}+\mu_2)!\\
\times [x_1^{\mu_1}] f_{k_1}(x_1)^{\nu}\,[y_2^{\mu_2}] f_{k_2}(y_2)^{N-\nu}\\
\times [y^{\mu_1+\mu_{1,2}}]f_{k_2}(y)^{\nu}\,[x^{\mu_{1,2}+\mu_2}]f_{k_1}(x)^{N-\nu}.
\end{multline}
So, by \eqref{Snu,mu=},
\begin{multline}\label{Snu,mu=expl}
S_{\nu,\boldsymbol\mu}=M!\binom{N}{\nu}\cdot\frac{(\mu_1+\mu_{1,2})!(\mu_{1,2}+\mu_2)!}
{\mu_{1,2}!}\\
\times [x_1^{\mu_1}] f_{k_1}(x_1)^{\nu}\,[y_2^{\mu_2}] f_{k_2}(y_2)^{N-\nu}\\
\times [y^{\mu_1+\mu_{1,2}}]f_{k_2}(y)^{\nu}\,[x^{\mu_{1,2}+\mu_2}]f_{k_1}(x)^{N-\nu}.
\end{multline}

To bound $S_{\nu,\boldsymbol\mu}$, we need an inequality
\begin{equation}\label{gammak}
[x^{a}] f_k(x)^b\le \frac{\gamma_k}{\sqrt{bx}}\,\frac{f_k(x)^b}{x^a},\quad
\forall\,x>0,
\end{equation}
where $\gamma_k$ depends on $k$ only. The proof follows from the Cauchy integral formula
\[
[x^{a}] f_k(x)^b=\frac{1}{2\pi i}\oint_{|z|=x}\frac{f_k(z)^b}{z^{a+1}}\,dz,
\]
and an inequality
\[
|f_k(xe^{i\theta})|\le f_k(x)\exp\left(x\,\frac{\cos\theta-1}{k+1}\right), \quad (x>0),
\]
see \cite{Pit1}. In addition, $f_k(x)$ is log-concave for $x>0$ since
\[
(\log f_k(x))'=\frac{f_{k-1}(x)}{f_k(x)}=1+\left(\sum_{j\ge k}\frac{(k-1)!}{j!}x^{j-k+1}\right)^{-1}
\]
decreases with $x$.

Using \eqref{gammak} and  log-concavity of $f_k$, we have
\begin{align*}
[x_1^{\mu_1}]f_{k_1}(x_1)^{\nu}\,[x^{\mu_{1,2}+\mu_2}]f_{k_1}(x)^{N-\nu}&\le\frac{\gamma_{k_1}^2}{\sqrt{\nu x_1(N-\nu)x}}
\frac{f_{k_1}(x_1)^{\nu}}{x_1^{\mu_1}}\,\frac{f_{k_1}(x)^{N-\nu}}{x^{\mu_{1,2}+\mu_2}}\\
&\le\frac{\gamma_{k_1}^2}{\sqrt{\nu x_1(N-\nu)x}}\frac{f_{k_1}\left(\frac{\nu}{N}x_1+\frac{N-\nu}{N}x\right)^N}{x_1^{\mu_1}\,x^{\mu_{1,2}+\mu_2}},
\end{align*}
for all $x_1>0,\,x>0$. By $\mu_1+\mu_{1,2}+\mu_2=M$ and \eqref {zi,zo;roots}, we easily obtain 
that the last fraction attains its minimum at
\[
x_1=\frac{Nz_i}{M}\cdot\frac{\mu_1}{\nu},\quad x=\frac{Nz_i}{M}\cdot\frac{\mu_{1,2}+\mu_2}{N-\nu},
\]
and the minimum itself is
\[
\frac{f_{k_1}(z_i)^N}{z_i^M}\,\left(\frac{M}{N}\right)^M\left(\frac{\nu}{\mu_1}\right)^{\mu_1}\left(\frac{N-\nu}{\mu_{1,2}+\mu_2}\right)^{\mu_{1,2}+\mu_2}.
\]
Therefore
\begin{equation}\label{easy}
\begin{aligned}
&[x_1^{\mu_1}]f_{k_1}(x_1)^{\nu}\,[x^{\mu_{1,2}+\mu_2}]f_{k_1}(x)^{N-\nu}\\
&\le \frac{\gamma_{k_1}^2}{\sqrt{\mu_1(\mu_{1,2}+\mu_2)}}\frac{f_{k_1}(z_i)^N}{z_i^{M+1}}\!\left(\frac{M}{N}\right)^{M+1}\!\!\left(\frac{\nu}{\mu_1}\right)^{\mu_1}\!\!\left(\frac{N-\nu}{\mu_{1,2}+\mu_2}\right)^{\mu_{1,2}+\mu_2}\!.
\end{aligned}
\end{equation}
Similarly
\begin{equation}\label{last}
\begin{aligned}
&[y_2^{\mu_2}] f_{k_2}(y_2)^{N-\nu} [y^{\mu_1+\mu_{1,2}}]f_{k_2}(y)^{\nu}\\
&\le\! \frac{\gamma_{k_2}^2}{\sqrt{\mu_2(\mu_{1,2}+\mu_1)}}\frac{f_{k_2}(z_o)^N}{z_o^{M+1}}\!\left(\frac{M}{N}\!\right)^{M+1}\!\!\!\left(\!\frac{N-\nu}{\mu_2}\!\right)^{\mu_2}
\!\!\left(\frac{\nu}{\mu_1+\mu_{1,2}}\!\right)^{\mu_1+\mu_{1,2}}\!\!.
\end{aligned}
\end{equation}
Combining \eqref{easy} and \eqref{last}, we get from \eqref{Snu,mu=expl} and 
$b!\le \text{const } b^{1/2}(b/e)^b$  that: 
\begin{equation}\label{Snu,mu=O}
\begin{aligned}
S_{\nu,\boldsymbol\mu}=&O\Biggl(\!\binom{N}{\nu}\!\left(\frac{M}{eN}\right)^{2M}\frac{f_{k_1}(z_i)^Nf_{k_2}(z_o)^N}{z_i^M\,z_o^M}\\
&\times\,\binom{M}{\boldsymbol\mu}\nu^{2\mu_1+\mu_{1,2}}(N-\nu)^{\mu_{1,2}+2\mu_2}\!\Biggr),
\end{aligned}
\end{equation}
were $\binom{M}{\boldsymbol\mu}$ stands for the trinomial coefficient.

Our next step is to add up the bounds \eqref{Snu,mu=O} for all $\boldsymbol\mu$ meeting
the constraints \eqref{mularge} and $\mu_1+\mu_{1,2}+\mu_2=M$. As it turns out, we will not lose anything by paying
attention to a single constraint $\mu_1\ge k_1\nu$ in \eqref{mularge}. Intuitively this is because the dominant
contribution to the sum comes from small $\nu$ and $\mu_1$, in which case this constraint is
most stringent among those in \eqref{mularge}.
So, using Chernoff's method, we introduce $u>1$ and bound
\begin{multline}\label{u1}
\sum_{\boldsymbol\mu\text{ meets }\eqref{mularge},\atop \mu_1+\mu_{1,2}+\mu_2=M}\binom{M}{\boldsymbol\mu}
\nu^{2\mu_1+\mu_{1,2}}(N-\nu)^{\mu_{1,2}+2\mu_2}\\
\le u^{-k_1\nu}
\sum_{\mu_1+\mu_{1,2}+\mu_2=M}\binom{M}{\boldsymbol\mu}
u^{\mu_1}\nu^{2\mu_1+\mu_{1,2}}(N-\nu)^{\mu_{1,2}+2\mu_2}\\
=u^{-k_1\nu}\bigl(\nu^2u+\nu(N-\nu)+(N-\nu)^2\bigr)^M\\
=\bigl(\nu^2+N(N-\nu)\bigr)^M.
\end{multline}
It follows then from \eqref{Snu,mu=O} and \eqref{mathcal SNM=} that
\begin{equation}\label{frac=O}
\frac{\sum_{\boldsymbol\mu}S_{\nu,\boldsymbol\mu}}{S_{N,M}}=
O\left(\binom{N}{\nu}u^{-k_1\nu}\bigl(1-\rho+u\rho^2\bigr)^M\right),\quad \rho:=\frac{\nu}{N}.
\end{equation}
Denoting $\sigma=M/N$, and using $\sigma>\max\{k_1,k_2\}$, the RHS of \eqref{frac=O} attains its absolute minimum at
\[
u_{\text{min}}=\frac{k_1(1-\rho)}{\rho(\sigma-k_1\rho)},
\]
and $u_{\text{min}}\ge 1$ iff 
\begin{equation}\label{u1>1}
\rho\le\rho^*:=\frac{2k_1}{k_1+\sigma+\sqrt{(\sigma-k_1)(\sigma+3k_1)}}.
\end{equation}
Recall that we consider $\rho=\nu/N\le 1/2$. All those $\rho$ will meet the constraint \eqref{u1>1}
iff $\rho^*\ge 1/2$, which is equivalent to $\sigma\le 3k_1/2$.

So, for  $\sigma\le 3k_1/2$,  from \eqref{mathcal SNM=} we get
\begin{equation}\label{???}
\begin{aligned}
\frac{\sum_{\boldsymbol\mu}S_{\nu,\boldsymbol\mu}}{S_{N,M}}
&=O\left[\binom{N}{\nu}\frac{\left(\tfrac{\sigma-\sigma\rho}{\sigma-k_1\rho}\right)^M}
{\left(\frac{k_1(1-\rho)}{\rho(\sigma-k_1\rho)}\right)^{k_1\nu}}\right]\\
&=O\bigl(\nu^{-1/2}\,e^{N H(\rho,\sigma)}\bigr),\\
H(\rho,\sigma)&:=\rho\log\frac{1}{\rho}+(1-\rho)\log\frac{1}{1-\rho}
+ \sigma\log\left(\tfrac{\sigma-\sigma\rho}{\sigma-k_1\rho}\right)\\
&\quad -k_1\rho\log \left(\frac{k_1(1-\rho)}{\rho(\sigma-k_1\rho)}\right).
\end{aligned}
\end{equation}
Since $\sigma>k_1$, we have
\[
\frac{\partial H(\rho,\sigma)}{\partial\sigma}=\log\frac{\sigma-\sigma\rho}{\sigma-k_1\rho}<0,\quad
\bigl(\rho\in (0,1]\bigr).
\]
So, for all $\rho\in (0,1/2]$, we have
\begin{equation}\label{??}
\begin{aligned}
H(\rho,\sigma)&< H(\rho,k_1)=(1-\rho)\log\frac{1}{1-\rho}-(k_1-1)\rho\log\frac{1}{\rho}\\
&\le (1-\rho)\log\frac{1}{1-\rho}-\rho\log\frac{1}{\rho}\le 0.
\end{aligned}
\end{equation}
for all $\rho\in (0,1/2]$, whence for $\rho\in \bigl(0, \min\{\rho^*,1/2\}\bigr)$.

Let $\sigma>3k_1/2$, so that $\rho^*<1/2$.  In this case the bound \eqref{???} continues to
hold for $\rho\in (0,\rho^*]$. However,  for $\rho\in [\rho^*,1/2]$ the RHS of  \eqref{u1}  attains its minimum at $u_{\text{min}}=1$,
and, instead of the bound \eqref{???}, we get 
\begin{equation}\label{????}
\begin{aligned}
\frac{\sum_{\boldsymbol\mu}S_{\nu,\boldsymbol\mu}}{S_{N,M}}&=O\left[\binom{N}{\nu}(1-\rho+\rho^2)^M\right]\\
&=O\bigl(\nu^{-1/2}e^{NK(\rho,\sigma)}\bigr),\\
K(\rho,\sigma)&:=\rho\log\frac{1}{\rho}+(1-\rho)\log\frac{1}{1-\rho}+\sigma\log(1-\rho+\rho^2).
\end{aligned}
\end{equation}
Observe that
\begin{equation}\label{double}
\frac{\partial^2 K(\rho,\sigma)}{\partial\rho^2}=-\frac{1}{\rho(1-\rho)}+\sigma\frac{1+2\rho(1-\rho)}
{\bigl(1-\rho(1-\rho)\bigr)^2}
\end{equation}
is increasing on $(0,1/2]$ as $\rho(1-\rho)$ is increasing. Therefore, as a function of $\rho$, $K(\rho,\sigma)$ is convex on $[\rho^*,1/2]$ iff $\tfrac{\partial^2 K(\rho,\sigma)}{\partial\rho^2}\ge 0$
at $\rho^*$, or equivalently, by \eqref{double}, iff 
\[
\rho^*(1-\rho^*)\ge \frac{2}{\sigma+2+\sqrt{\sigma^2+12\sigma}}.
\]
From the definition of $\rho^*$ in \eqref{u1} it follows that, as a function of $k_1$, $\rho^*$ is
increasing as long as $\rho^*\le 1/2$. So the condition above holds for all $k_1\in [2, 2\sigma/3]$, if it does for $k_1=2$, in which case $\sigma\ge 3$. An elementary algebraic verification
does the job. 

Now $\tfrac{\partial K(\rho,\sigma)}{\partial\rho}=0$ at $\rho=1/2$;  so $K(\rho,\sigma)$, the convex function on $[\rho^*,1/2]$, attains its minimum at $\rho=1/2$. Therefore $K(\rho,\sigma)$ attains its maximum at the other end $\rho^*$, and by definition of $\rho^*$, 
\[
K(\rho^*,\sigma)=H(\rho^*,\sigma)<0.
\]
Combining \eqref{???}, \eqref{??}, \eqref{????} and the last inequality, we conclude
that $S_{N,M}^{-1}\sum_{\boldsymbol\mu}S_{\nu,\boldsymbol\mu}$ is uniformly exponentially small for all $\nu/N=\rho\in [\eps,1/2]$, $\eps>0$ being arbitrarily small, and for $\rho<\eps$
\begin{equation*}
\frac{\sum_{\boldsymbol\mu}S_{\nu,\boldsymbol\mu}}{S_{N,M}}
=O\bigl[\nu^{-1/2}\exp(-N(k_1-1)(\rho\log1/\rho +O(\rho)))\bigr].
\end{equation*}
So
\begin{equation}\label{sum1,N/2}
\sum_{\nu=2}^{N/2}\frac{\sum_{\boldsymbol\mu}S_{\nu,\boldsymbol\mu}}{S_{N,M}}=O(N^{-2(k_1-1)}).
\end{equation}
Likewise the expected number of sink sets of size in $[2, N/2]$ is $O(N^{-2(k_2-1)})$. We conclude that the uniformly random $(k_1,k_2)$-core is strongly connected with probability 
$1-O\bigl(N^{-2(k-1)}\bigr)$, $k:=\min\{k_1,k_2\}$.

Let us show that in fact the random core is $k$-strongly connected with probability 
$1-O\bigl(N^{-(k-1)}\bigr)$. That is, we want to show that for $t\in [1,k-1]$ w.h.p. there does not exist a 
partition $[N]=A_1\uplus A_2\uplus A_3$, with $|A_1|=\nu_1\in [2,(N-t)/2]$, $|A_2|=\nu_2$, and $|A_3|=\nu_3=t$ such that deletion of $A_3$ results in a digraph where $A_1$ is a source (sink) set. To do so we need to prove that the total number of sequences inducing such a partition of $[N]$ is
$o(S_{N,M})$ as $N\to\infty$. The argument is a natural extension of the proof of strong
connectedness. So we will focus on the new details.

Introduce $\mu_j$, $1\le j\le 3$,
the (generic) numbers of edges in the vertex sets $A_j$, and $\mu_{i,j}$, the number
of edges from the $A_i$ to the $A_j$. For $A_1$ to be 
a source set upon deletion of $A_3$ we must have $\mu_{2,1}=0$, and 
\[
\mu_1+\mu_{3,1}\ge k_1\nu_1,
\]
besides all other constraints, similar to those in \eqref{mularge}. Define $\mu_{i,i}=\mu_i$, and
\[
\mu_{\boldsymbol\cdot,j}=\sum_{i}\mu_{i,j},\quad \mu_{j,\boldsymbol\cdot}=\sum_{i}\mu_{j,i},
\]
Let $S_{\boldsymbol\nu,\boldsymbol\mu}$ be the total number of the sequences with these 
parameters. Analogously to \eqref{Snu,mu=},
\begin{equation}\label{Sboldnu,boldmu}
S_{\boldsymbol\nu,\boldsymbol\mu}=\binom{N}{\boldsymbol\nu}\binom{M}{\boldsymbol\mu}
\mathcal{S}_{\boldsymbol\nu,\boldsymbol\mu},\quad \binom{N}{\boldsymbol\nu}=\frac{N!}{\nu_1!
\nu_2!\nu_3!}.
\end{equation}
Here $\mathcal{S}_{\boldsymbol\nu,\boldsymbol\mu}$ is defined like $\mathcal{S}_{\nu,\boldsymbol\mu}$, with $[N]$ partitioned in three consecutive blocks $A_1$, $A_2$, $A_3$ of length $\nu_1$, $\nu_2$, and $\nu_3=t$.  Analogously to \eqref{long}, we have
\[
\mathcal S_{\boldsymbol\nu,\boldsymbol\mu}=\prod_{j=1}^3\mu_{j,\boldsymbol\cdot}!\,\mu_{\boldsymbol\cdot,j}!
\times \bigl[\xi_j^{\mu_{\boldsymbol\cdot,j}}\bigr] f_{k_1}(\xi_j)^{\nu_j}\,\bigl[\eta_j^{\mu_{j,\boldsymbol\cdot}}\bigr]f_{k_2}(\eta_j)^{\nu_j}.
\]
Here, like \eqref{easy}-\eqref{last},
\begin{align*}
&\prod_{j=1}^3\bigl[\xi_j^{\mu_{\boldsymbol\cdot,j}}\bigr] f_{k_1}(\xi_j)^{\nu_j}\le \gamma_{k_1}^3\frac{f_{k_1}(z_i)^N}{z_i^{M+3/2}}
\left(\frac{M}{N}\right)^{M+3/2} \prod_{j=1}^3\mu_{\boldsymbol\cdot,j}^{-1/2}\left(\frac{\nu_j}{\mu_{\boldsymbol\cdot,j}}\right)^{\mu_{\boldsymbol\cdot,
j}},\\
&\prod_{j=1}^3\bigl[\eta_j^{\mu_{j,\boldsymbol\cdot}}\bigr] f_{k_2}(\eta_j)^{\nu_j}\le \gamma_{k_2}^3\frac{f_{k_2}(z_o)^N}{z_o^{M+3/2}}\left(\frac{M}{N}\right)^{M+3/2}\prod_{j=1}^3\mu_{j,\boldsymbol\cdot}^{-1/2}\left(\frac{\nu_j}{\mu_{j,\boldsymbol\cdot}}\right)^{\mu_{j,\boldsymbol\cdot}}.
\end{align*}
Consequently
\[
\mathcal S_{\boldsymbol\nu,\boldsymbol\mu}=O\left(\frac{f_{k_1}(z_i)^Nf_{k_2}(z_o)^N}{z_i^M z_o^M}\left(\frac{M}{eN}\right)^{2M}
\prod_{j=1}^3\nu_j^{\mu_{j,\boldsymbol\cdot}+\mu_{\boldsymbol\cdot,j}}\right).
\]
This bound, combined with \eqref{Sboldnu,boldmu}, \eqref{mathcal SNM=} and $\sum_j(\mu_{j,\boldsymbol\cdot}+\mu_{\boldsymbol\cdot,j})=2M$, delivers
\begin{equation}\label{del}
\frac{S_{\boldsymbol\nu,\boldsymbol\mu}}{S_{N,M}}
=O\left(\binom{N}{\boldsymbol\nu}\binom{M}{\boldsymbol\mu}
\prod_{j=1}^3\left(\frac{\nu_j}{N}\right)^{\mu_{j,\boldsymbol\cdot}
+\mu_{\boldsymbol\cdot,j}}\right).
\end{equation}
Next, analogously to the case $\nu_3=t=0$, we have: for $u\ge 1$ {\it and\/}  $u=O(N/\nu_1)$,
\begin{align*}
&\sum_{\|\boldsymbol\mu\|=M,\atop \mu_{2,1}=0,\,
\mu_1+\mu_{3,1}\ge k_1\nu_1}\!\!\!\binom{M}{\boldsymbol\mu}\,
\prod_{j=1}^3\left(\frac{\nu_j}{N}\right)^{\mu_{j,\boldsymbol\cdot}+\mu_{\boldsymbol\cdot,j}}\\
&\le u^{-k_1\nu_1}\sum_{\|\boldsymbol\mu\|=M}u^{\mu_1+\mu_{3,1}}
\binom{M}{\boldsymbol\mu}\,\prod_{j=1}^3\left(\frac{\nu_j}{N}\right)^{\mu_{j,\boldsymbol\cdot}+\mu_{\boldsymbol\cdot,j}}\\
&=O\left[u^{-k_1\nu_1}\left(u\frac{\nu_1(\nu_1+\nu_3)}{N^2}+\frac{\nu_2+\nu_3}{N}\right)^M
\right]\\
&=O\left(u^{-k_1\nu_1}(1-\rho+u\rho^2)^M\right),\quad \rho:=\frac{\nu_1}{N}.
\end{align*}
So
\begin{equation}\label{3frac}
\begin{aligned}
\frac{\sum_{\boldsymbol\mu}S_{\boldsymbol\nu,\boldsymbol\mu}}{S_{N,M}}&=O\left(
\binom{N}{\boldsymbol\nu}u^{-k_1\nu_1}(1-\rho+u\rho^2)^M\right)\\
&=O\left(N^t\binom{N-t}{\nu_1}u^{-k_1\nu_1}(1-\rho+u\rho^2)^M\right).
\end{aligned}
\end{equation}
From this moment on we reason almost like in the argument following \eqref{frac=O}, and end up
with
\[
\sum_{\nu=2}^{(N-t)/2}\min_{u\ge 1}\frac{\sum_{\boldsymbol\mu}S_{\boldsymbol\nu,\boldsymbol\mu}}{S_{N,M}}=O\bigl(N^t\cdot N^{-2(k_1-1)}\bigr)=O\bigl(N^{-(k_1-1)}\bigr).
\]
So, for $t<k_1$ ($t<k_2$ resp.), with probability $1-O\bigl(N^{-(k_1-1)}\bigr)$ 
($1-O\bigl(N^{-(k_2-1)}\bigr)$ resp.)
there is no set $A$ of cardinality in $[2,(N-t)/2]$ which becomes a source set (sink set resp.) upon deletion of $t$ vertices from $[N]\setminus A$.

The proof of Theorem \ref{strcon} is complete.
\bigskip

{\bf Acknowledgment.\/}  I am grateful to Dan Poole for formulating succinctly his thought-provoking conjecture and stopping me from pursuing false leads.

\end{document}